\documentclass[11 pt,fullpage]{article}
\usepackage{hyperref}
\usepackage{ graphicx, mathrsfs}
\usepackage{amssymb,amsmath,amsthm}
\usepackage{mathrsfs,dsfont}

\theoremstyle{plain}
\newtheorem{theorem}{Theorem}

\newtheorem{corollary}{Corollary}

\newtheorem*{remark}{Remark}
\numberwithin{equation}{section}
\numberwithin{equation}{section}

\renewcommand\o{\omega}

\newcommand{\Real}{\mathbb R}
\newcommand{\T}{\mathbb T}
\newcommand{\Integer}{\mathbb Z}
\newcommand{\norm}[1]{\|#1\|}
\newcommand{\abs}[1]{\left\vert#1\right\vert}
\newcommand{\set}[1]{\left\{#1\right\}}

\newcommand{\grad}{\nabla}

\newcommand{\jap}[1]{\langle #1 \rangle}

\begin{document}

\title{Asymptotic stability for the Couette flow in the 2D Euler equations} 

\author{Jacob Bedrossian\footnote{\textit{jacob@cims.nyu.edu}, Courant Institute of Mathematical Sciences. Partially supported by NSF Postdoctoral Fellowship in Mathematical Sciences, DMS-1103765} \, and Nader Masmoudi\footnote{\textit{masmoudi@cims.nyu.edu}, Courant Institute of Mathematical Sciences. Partially supported by NSF  grant}}

\date{\today}
\maketitle

\begin{abstract}
In this expository note we discuss our recent work \cite{BM13}  on the nonlinear 
asymptotic  stability of shear flows in the 2D Euler equations of ideal, incompressible flow. 
In that work it is proved that perturbations to the Couette flow which are small in a suitable regularity class 
 converge strongly in $L^2$ to a shear flow which is close to the Couette flow. 
Enstrophy is mixed to small scales by an almost linear evolution and is generally lost in the weak limit as $t \rightarrow \pm\infty$. 
In this note we discuss the most important physical and mathematical aspects of the result and the key ideas of the proof. 
\end{abstract} 

\section{Introduction}

We consider the 2D Euler system in the vorticity formulation with a background shear flow: 
\begin{equation} \label{def:2DEuler}
\left\{
\begin{array}{l}
  \omega_t + y\partial_x\omega + U \cdot \grad \omega = 0, \\ 
  U  = \grad^{\perp}(\Delta)^{-1} \omega,  \quad \quad \omega(t=0) =\omega_{in}. 
\end{array}
\right. 
\end{equation}
Here,  $(x,y) \in \mathbb T \times \Real$, $\grad^\perp = (-\partial_y,\partial_x)$ and $(U,\omega)$ are periodic in the $x$ variable with period normalized to $2\pi$.   
The physical velocity is $(y,0) + U$ where $U = (U^x,U^y)$ denotes the velocity perturbation and the total vorticity is $-1 + \omega$. 
We denote the streamfunction $\psi = \Delta^{-1}\omega$. 

The field of hydrodynamic stability has a long history starting in the nineteenth century. 
One of the oldest problems considered is the stability and instability of shear flows, dating back to, for example, Rayleigh \cite{Rayleigh80} and Kelvin \cite{Kelvin87}. 
For the case considered here, the solution to the linearization of \eqref{def:2DEuler} can be found in the work of Kelvin \cite{Kelvin87} although the solution in the inviscid case was not specifically analyzed until the work of Orr in 1907 \cite{Orr07}.
See \cite{BM13} for a discussion of the history of \eqref{def:2DEuler} and its relationship with the wider field of hydrodynamic stability (as well as many related references).  
We first state the result of \cite{BM13} and then attempt to elucidate some of the interesting physical and mathematical concepts which are involved in the proof.  
The relationship with Landau damping in the Vlasov equations of plasma physics and the recent work of Mouhot and Villani \cite{MouhotVillani11} is discussed as well. 

\section{Asymptotic Stability} 
We are interested in the long time behavior of  \eqref{def:2DEuler} for a  small initial 
 perturbation $ \omega_{in}$. 
In particular, we are interested in studying the \emph{asymptotic stability} of shear flows, that is
 showing that all sufficiently small perturbations in a suitable regularity class converge to a shear flow in the sense: $y + U(t,x,y) \rightarrow U_\infty(x,y) = (y + u_\infty(y),0)$ as $t \rightarrow \infty$. 
The initial vorticity is taken in a Gevrey space of class $1/s$ for $s> 1/2$ \cite{Gevrey18}; although this regularity class is slightly unusual, this restriction on the initial data arises naturally from the weakly nonlinear effects.
See \cite{BM13} for some speculation on whether this requirement is sharp (it is shown in \cite{LinZeng11} that the regularity requirement must be at least $H^{3/2}$ but this is a huge gap; this work is discussed more below in \S\ref{sec:Nonlinear}). 
We note that the analogous space for the Vlasov equations with Coloumb/Newton interaction is Gevery-3 (e.g. $s = 1/3$) \cite{MouhotVillani11}. 
 Our main result is
\begin{theorem} \label{thm:Main} 
For all $1/2 < s \leq 1$, $\lambda_0 > \lambda^\prime > 0$ there exists an $\epsilon_0  = \epsilon_0(\lambda_0,\lambda^\prime,s) \leq 1/2$ such that for all $\epsilon  \leq  \epsilon_0   $  
 if $\omega_{in}$ satisfies 
  $\int \omega_{in} dx dy= 0$, $\int \abs{y\omega_{in}(x,y)} dx dy < \epsilon$ and   
\begin{align*}  
\norm{\omega_{in}}^2_{\lambda_0} = \sum_k\int \abs{\hat{\omega}_{in}(k,\eta)}^2 e^{2\lambda_0 \abs{k,\eta}^{s}} d\eta \leq  \epsilon^2,
\end{align*}  
then there exists $f_\infty$ with $\int f_\infty dxdy = 0$ 
 and $\norm{f_\infty}_{\lambda^\prime} \lesssim \epsilon$ such that 
\begin{equation} \label{main-omega} 
\norm{\omega(t,x + ty + \Phi(t,y),y) - f_\infty(x,y) }_{\lambda^\prime} \lesssim \frac{\epsilon^2}{\jap{t}}, 
\end{equation} 
where $\Phi(t,y)$ is given explicitly by 
\begin{align} 
\Phi(t,y) = \frac{1}{2\pi}\int_0^t\int_\T U^x(s,x,y) dx ds = u_\infty(y)t + \theta(t,y),    \label{def:phi}
\end{align} 
with $u_\infty = \partial_y \partial_{yy}^{-1}\frac{1}{2\pi}\int_\T f_\infty(x,y) dx$ and 
  $| \theta(t,y) | \lesssim   \epsilon^2\abs{\log t}$.
Moreover, the velocity field $U$ decays as
\begin{subequations} \label{ineq:damping}
\begin{align} 
\norm{\frac{1}{2\pi}\int U^x(t,x,\cdot) dx - u_\infty}_{\lambda^{\prime}} &\lesssim \frac{\epsilon^2}{\jap{t}}, \label{ineq:xdamping_slow} \\ 
\norm{U^x(t) - \frac{1}{2\pi}\int U^x(t,x,\cdot) dx}_{L^2} &\lesssim \frac{\epsilon}{\jap{t}}, \label{ineq:xdamping} \\ 
\norm{U^y(t)}_{L^2} & \lesssim \frac{\epsilon}{\jap{t}^2}. \label{ineq:ydamping}
\end{align}
\end{subequations}
By time-reversibility, statements analogous to \eqref{def:phi} and \eqref{ineq:damping} also hold backward in time for some $f_{-\infty}$ (which generally has no reason to be equal to $f_{\infty}$). 
\end{theorem}

\begin{remark} 
As in the scattering of nonlinear dispersive equations and Landau damping in the Vlasov equations \cite{MouhotVillani11}, there is no simple formula to determine $f_\infty$, which is chosen by the nonlinear evolution. However, it follows from the proof that $\norm{\omega_{in} - f_\infty}_{\lambda^\prime} \lesssim \epsilon^2$. %as the nonlinear effects are one order weaker than the linear effects (in amplitude).  
\end{remark}

The `inviscid damping' \eqref{ineq:damping} expresses the strong convergence of the primary observable, namely the velocity field, to a shear flow.  
On the other hand, \eqref{main-omega} describes the asymptotic evolution of the vorticity as a logarithmic correction to passive transport in a shear flow: 
\begin{align*}
\omega(t,x,y) \sim f_\infty(x - ty - tu_\infty(y) - \theta(t,y),y), 
\end{align*}
and hence is expressing \emph{weak} convergence of the vorticity and the transfer of enstrophy to small scales, some of which may be lost as $t \rightarrow \pm \infty$. 
In fact there is the following corollary which follows easily from the proof of Theorem \ref{thm:Main} (see \cite{BM13}).  
\begin{corollary} \label{cor1}
There exists an open set of smooth solutions to \eqref{def:2DEuler} for which $\set{\omega(t)}_{t \in \Real}$ is not pre-compact in $L^2$ as $t \rightarrow \pm\infty$. In particular, $\omega(t) \rightharpoonup \omega_\infty = \frac{1}{2\pi}\int_\T f_\infty(x,y) dx$ and in general $\norm{\omega_\infty}_2 < \norm{\omega(t)}_2$. 
\end{corollary} 
This loss of enstrophy to high frequencies was widely expected but to our knowledge, our work  
 \cite{BM13}  is the first rigorous confirmation of the effect in the 2D Euler equations (albeit in a somewhat specific setting). 
Indeed, one of the primary motivations for our work was to rigorously study this effect on the nonlinear level, which is thought to be a fundamental mechanism in the 2D Euler equations connected to the meta-stability of coherent structures and 2D turbulence \cite{Kraichnan67,Gilbert88,Shnirelman12,GSV13}. 
%Another related corollary to Theorem \ref{thm:Main} is that it provides an open set of smooth solutions for which $\norm{\omega(t)}_{H^k} = O(t^k)$. 

In his 1907 paper, Orr \cite{Orr07} studied the linearization of 2D Euler around the planar Couette flow (among several other configurations), seeking to reconcile the linear stability with the instabilities often seen in experiments (see \cite{LiLin11} and the references therein for more discussion on this disagreement, sometimes called the `Sommerfeld paradox').
Orr's observations play a major role in our work and can be summarized in modern terminology (and adapted to our infinite-in-$y$ setting) as follows.  
Given a disturbance in the vorticity, on the linear level it is simply advected by the background shear flow: $\omega(t,x,y) = \omega_{in}(x-ty,y)$. 
If one changes coordinates to $z = x-ty$ then the stream-function $\phi(t,z,y)$ in these variables solves $\partial_{zz}\phi +
  (\partial_y - t\partial_z)^2\phi = \omega_{in}$. 
On the Fourier side, $(z,y) \to (k,\eta) \in \Integer \times \Real$,
\begin{align} \label{orr-cri} 
\hat{\phi}(t,k,\eta) = -\frac{\hat{\omega}_{in}(k,\eta)}{k^2 + \abs{\eta - kt}^2}. 
\end{align}  
From \eqref{orr-cri}, Orr made two important observations, together known now as the \emph{Orr mechanism}.  
First, he formally identified the decay rates \eqref{ineq:xdamping}\eqref{ineq:ydamping}. 
His work pre-dated Sobolev spaces, but we would now see this observation as the fundamental inequality 
\begin{align}
\norm{P_{\neq 0}\phi}_{H^N} \lesssim \frac{\norm{\omega_{in}}_{H^{N+2}}}{1 + t^2}, \label{ineq:phidecay}
\end{align}
where $H^N$ denotes the Sobolev space of order $N$ and $P_{\neq 0}\phi = \phi - \frac{1}{2\pi}\int \phi dz$, the projection onto non-zero frequencies in $z$ (the loss of the power of $t$ in going to \eqref{ineq:xdamping} is due to the time-dependence of the change of variables). 
Second, Orr identified the possibility for a large transient growth in the kinetic energy of the disturbance. 
Specifically, this transient growth occurs forward in time for those modes which satisfy $t_c = \frac{\eta}{k} > 0$ (referred to by Orr as `critical times').
Physically, these correspond to waves tilted against the shear which are being deformed to larger wave-lengths. 
Orr suggested this latter effect as an explanation for the \emph{practical instability} observed in experiments.
Indeed, it is seen from \eqref{orr-cri} that the linearized problem is \emph{unstable} with respect to the kinetic energy of the perturbation.
Moreover, upon consideration we see that \eqref{ineq:phidecay} is essentially optimal in the sense that we cannot expect decay without the loss of two derivatives (and we cannot, in general, get faster decay by paying more). 
This highlights an important theme underlying our work and that of \cite{MouhotVillani11}, which is that \emph{decay costs regularity}. 
Indeed, that we must pay something is not surprising when we reflect on the fact that \eqref{def:2DEuler} can be seen as a Hamiltonian system and hence we should only expect asymptotic stability in a norm weaker than that taken on the initial data. 
This also brings up another interesting property discussed further below, which is that the behavior in Theorem \ref{thm:Main} and in Landau damping is intrinsically \emph{infinite dimensional} and has no direct analogue in finite dimensional systems. 
 
The decay predicted by the Orr mechanism is due to the transfer of enstrophy to small scales (which yields decay of the velocity through the Biot-Savart law) 
and the transient growth can be understood as the time-reversed phenomenon: the transfer of enstrophy from small scales to large scales (see \cite{BM13, Boyd83,Lindzen88} for further discussion).
This transfer of information to small scales is rightly considered a hydrodynamic analogue of Landau damping in the collisionless Vlasov equations of plasma physics, which is the origin of the term \emph{inviscid damping} \cite{BouchetMorita10,SchecterEtAl00,Briggs70,BM95}.
The nonlinear theory of Landau damping was greatly advanced only recently by the work of Mouhot and Villani \cite{MouhotVillani11} (see also \cite{CagliotiMaffei98,HwangVelazquez09}), who showed that the Vlasov equations undergo Landau damping similar to that predicted on the linear level by Landau in 1946 \cite{Landau46}. 
The mathematical relationship between our work and \cite{MouhotVillani11} is discussed in more detail in \cite{BM13}.  
However, let us here point out several important differences between Landau damping in its simplest setting and inviscid damping for \eqref{def:2DEuler}. 
During inviscid damping, the velocity field in \eqref{def:2DEuler} does not converge back to Couette flow, but instead only converges to some nearby shear flow, whereas the electric field in the Vlasov equations converges rapidly to zero during Landau damping.
This is responsible for the fact that the final almost-linear evolution in \eqref{main-omega} depends on the solution itself.
This kind of `quasi-linearity' proves to be a very serious difficulty in studying inviscid damping on the nonlinear level.
Another key difference is that the decay \eqref{ineq:damping} is not even integrable for the $x$ component of the velocity whereas in the Vlasov equations, the decay is exponentially fast for analytic perturbations.
A third notable difference is the fact that the non-local law which gives rise to velocity in \eqref{def:2DEuler} does not distinguish between the $x$ and $y$ variables; that structure is imposed by the background flow. However, for the Vlasov equations, the electric field only depends on the density, which is a function of $x$ alone.  

Despite the differences, both Landau damping and inviscid damping are examples of the same general phenomenon observed in a number of infinite-dimensional Hamiltonian systems, now sometimes referred to as \emph{phase-mixing}, or \emph{continuum damping} (in reference to the spectra of the associated linearized operators). 
Indeed such behavior arises in numerous fluid mechanics applications \cite{BMSEI95,BouchetMorita10,CV13,Boyd83,Gilbert88,SchecterEtAl00,YuDriscoll02,YuDriscollONeil}, ideal MHD \cite{TataronisGrossman73}, quantum mechanics \cite{HagstromMorrison} and even biology \cite{StrogatzEtAl92}; see also the review article \cite{BMT13} and the references therein.  
On the linear level, the effect is directly tied to the continuous spectrum in the linear operator. 
Normal form-like transforms very similar to the generalized eigenfunctions employed by van Kampen \cite{VKampen55} show that one can re-write a variety of phase-mixing linear problems as a continuum of decoupled harmonic oscillators \cite{BM95,BMSEI95,Morrison98,Morrison00,BM02,HagstromMorrison,CV13}.
These decompositions essentially provide a natural infinite-dimensional analogue of action-angle coordinates for completely integrable Hamiltonian systems in finite dimensions.  
%This prompts some authors to draw a connection between nonlinear phase-mixing results and KAM theory (for example \cite{MouhotVillani11} \textbf{[Are there others?]}) but the precise nature of this connection is still far from clear. 
Phase mixing can also be seen from the RAGE theorem \cite{cycon1987,ReedSimonIII}, which implies the decay of compact operators applied to the linear evolution if the spectrum is purely absolutely continuous. 
Indeed, while the RAGE theorem originated in quantum mechanics, it has already seen applications to related mixing phenomena in fluids \cite{ConstantinEtAl08}. 

Aside from the RAGE theorem, several of the properties discussed in the previous paragraph are also shared by dispersive/wave phenomena, as already pointed out by several authors, for example \cite{Degond86,CagliotiMaffei98,MouhotVillani11}. 
Dispersion is also a distinctly infinite-dimensional effect commonly seen in Hamiltonian systems and like phase-mixing, is directly connected to the continuous spectrum of the linear operator.
Moreover, in both phase-mixing and scattering the long-time behavior is governed by the linear operator (or a modified version due to `long-range' effects as in Theorem \ref{thm:Main}; see e.g. \cite{GV00,Nakanishi02,IonescuPusateri13} and the references therein)
however the exact evolution cannot generally be characterized by simple principles since (by time-reversibility) it retains all of the information of the nonlinear dynamics. 
By time-reversibility, in either we may equivalently study forward and backward limits $t$ goes to $\pm \infty$.  
Also in both cases, there is a transient growth which can be seen if we consider the evolution between a time close to $-\infty$ and time $0$ (for example). 
However, one should not be tempted to see mixing as dispersion: the physical mechanisms are very different.
In dispersive systems, decay (for example in $L^\infty$) arises due to the fact that wave packets travel at different group velocities, which tends to spread the information out to spatial infinity.
Such a transfer of information typically costs \emph{spatial localization} while Sobolev norms are conserved.  %(for example, the $L^\infty$ decay of $L^1$ data for the free Schr\"odinger equation).   
However, mixing transfers information to infinity \emph{in frequency}, often preserving $L^p$ norms while causing observables obtained by compact operators to decay (provided there is regularity to pay).

\section{Physical and Mathematical aspects of the proof} 
We will only discuss the main new mathematical ideas in the proof and the associated physical intuition. 
The details and a more in-depth discussion are given in \cite{BM13}.
There are several key ingredients and technical difficulties:
\begin{enumerate} 
\item A change of variables that adapts to the solution as it evolves. 
Unlike the ideas of the  ``profile'' used in dispersive wave equations (see e.g. \cite{GMS12}) and of the ``gliding regularity'' in \cite{MouhotVillani11}, here it is crucial that the coordinate transformation depends on the solution.  
This seems necessary due to the `quasi-linearity' expressed in \eqref{main-omega}.  
\item The design of a special norm which loses regularity in a very precise way adapted to the Orr critical times and the associated weakly nonlinear effects. The norm is built using a reduced ``toy model'', which estimates the worst possible cascade of information to high frequencies (in the transformed variables).
Here it is estimated that Gevrey-2 regularity is potentially lost due to this cascade as $t \rightarrow \infty$, 
and hence here the restriction $s \geq 1/2$ arises.

\item A nonlinear energy estimate which uses the special norm to control the growth predicted by the toy model as well as subtle quasilinear effects that arise due to the fact that the coordinate transformation depends on the solution. 
This is expressed primarily by the fact that the coefficients of the Laplacian in the new variables, denoted $\Delta_t$ below, depend on a time-averaged projection of the vorticity.  
An additional energy estimate is also required to prove the convergence of the background shear flow, expressed in the time derivative of the coordinate transform. 
These latter complications currently have no analogue in the theory of Landau damping.

\item A new elliptic estimate for a degenerate, semilinear operator which represents the Laplacian expressed in the new coordinate system.
A main technical part of our proof is to gain regularity or decay from inverting it, despite the loss of ellipticity at the critical times as in \eqref{orr-cri}.  
The norm we design is (not coincidentally) well-adapted to exactly compensating for this loss of ellipticity which arises as a small denominator in the toy model (see \eqref{toy} below).   

 \end{enumerate} 

\subsection{Coordinate transformation}\label{chv}
As discussed above, from \eqref{main-omega}, we see two primary difficulties: $\o$ is not quite asymptotic to a linear evolution due to the logarithmic phase correction (analogous to modified scattering in dispersive equations \cite{GV00,Nakanishi02,IonescuPusateri13}) and more crucially, the shear flow, and hence linear evolution, itself is determined by the solution, and so can be regarded as an analogue of `quasilinear scattering'.
This is especially troublesome since as in the Vlasov equations, one can only get uniform estimates in norms adapted to the linear evolution \cite{MouhotVillani11}. 
Since we prescribe data at $t = 0$, we have no idea what $u_\infty$ or the logarithmic correction are; our way of dealing with this 
 lack of information  is to choose a coordinate system which adapts to the solution 
 and converges to the expected form as $t \rightarrow \infty$. 
The coordinates we use are
\begin{subequations} \label{def:zv}
\begin{align} 
z(t,x,y) & = x - tv \\ 
v(t,y)   & = y + \frac{1}{t}\int_0^t < U^x >(s,y) ds, \label{def:v}
\end{align}      
\end{subequations}
where $ < U^x >$ denotes the $x$-average of $U^x$.
The variable $v$ is named as such since it is an approximation of the shear velocity. 
Hence, the shift in $z$ is a `semi-Lagrangian' coordinate which approximately follows the background shear flow.  
This kind of shift is classical for studying Couette flow, dating back at least to the paper of Orr \cite{Orr07}, although \eqref{def:zv} is more nonlinear. 
The particular choice  $y \to v$ is to ensure that the Biot-Savart law is in a form amenable to Fourier analysis in the variables $(z,v)$. 

Define $f(t,z,v) = \omega(t,x,y)$ and the transformed stream-function $\phi(t,z,v) = (\Delta^{-1}_{x,y}\omega)(t,x,y)$. Then one derives (see \cite{BM13}),
\begin{align} \label{Euler2}
 \left\{
\begin{array}{l} 
\partial_tf + u \cdot \grad f = 0,  \\
u =  (0, [\partial_t v] )   + v^\prime \grad^\perp P_{\neq 0} \phi,  \\ 
\phi =\Delta_t^{-1}[f], 
\end{array}
\right.
\end{align} 
where we are denoting $[\partial_tv](t,v) = \partial_tv(t,y)$, $v^{\prime}(t,v) = \partial_y v(t,y)$ and $v^{\prime\prime}(t,v) = \partial_{yy} v(t,y)$. 
Actually, these functions admit relatively nice, direct representations in terms of $f$ (see \cite{BM13}).  
The Biot-Savart law is transformed into: 
\begin{align}
f & = \partial_{zz}\phi  + (v^\prime)^2\left( \partial_v - t \partial_z \right)^2\phi + v^{\prime\prime}\left(\partial_{v} - t \partial_z\right) \phi = \Delta_t \phi. \label{def:Deltat}
\end{align} 
The fact that the coefficients now depend on the (transformed) vorticity $f$ is the primary expression of the quasi-linearity 
and represents one of the main technical hurdles in the proof of Theorem \ref{thm:Main}. 

Given a smooth solution to \eqref{Euler2} such that  for all $t,v$, we have 
 $ \abs{\sup_{t}\frac{1}{t}\int_0^t f_0(s,v) ds}  < 1$,  we can recover a solution 
to the original system  \eqref{def:2DEuler} by an inverse function theorem.
Notice that $u$ is not divergence free and the dependence of $\phi$ on $f$ through $\Delta_t$ is significantly more subtle than in the original variables. 
The main advantage of \eqref{Euler2} is that $u$ formally has an integrable decay, indeed, if the solution is small and one is willing to pay four derivatives, the decay rate is formally  $O(t^{-2}\log t)$. 
If one has uniform control on the regularity of solutions to \eqref{Euler2}, Theorem \ref{thm:Main} 
can be recovered by an inverse function theorem argument in Gevrey spaces; see \cite{BM13} for how this is done. 

\subsection{Weakly nonlinear effects and the toy model} \label{sec:Nonlinear}
Since Orr's work, the unresolved fundamental question about the Couette flow is whether the Orr mechanism drives instability in the nonlinear 2D Euler equations or whether or not damping can still hold under some hypotheses. 
Lin and Zeng \cite{LinZeng11} showed that Theorem \ref{thm:Main} cannot hold if the initial vorticity is only $H^{s}$ for $s < 3/2$ by showing there exists a `cat's eye' steady state arbitrarily close to the Couette flow in sufficiently low regularities. They proved also the same result for Vlasov using the analogous BGK modes \cite{LZ11b}.
This work is telling in two ways: firstly, it shows that the nonlinear effect is closely tied to the regularity and secondly it hints at the expected dominant instability, which is the commonly observed vortex roll-up at the so-called `critical layer', which is where the disturbance `wave-speed' matches that of the background flow.
The fact that the proofs use the same instability is not a coincidence: it is possible to derive a `normal form' roughly analogous to the Landau equation, called the single-wave model, which formally describes this roll-up near the threshold of stability across several phase-mixing systems, including Vlasov, 2D shear flows, 2D vortices and the XY model of condensed matter physics; see the review article \cite{BMT13}. 

In our work, we are interested in showing that the evolution stays essentially linear uniformly in time and hence we need a very good understanding of the \emph{weakly nonlinear} effects.
It is a classical idea that transient growth in a linear problem can interact badly with the nonlinearity to trigger instability, for instance see the discussion in \cite{TTRD93}. The basic mechanism is as follows. Heuristically, in the weakly nonlinear regime we can imagine the solution as an interacting superposition of waves undergoing linear shear. 
Through the nonlinear term, each mode has a strong effect at its critical time during which it strongly forces the others, potentially putting information into modes which have not yet reached their critical time and are hence still growing. 
At a later time, these modes have a large effect and continue to excite other growing modes and so forth, perpetuating a so-called self-sustaining `nonlinear bootstrap' (see  \cite{TTRD93,BaggettEtAl,VMW98,Vanneste02} and the references therein for discussions in the fluid mechanics context). 
Since the measurable effect of a nonlinear interaction can occur long after the event, this mechanism permits the creation of controlled \emph{echoes}, in which the electric field of the plasma, or kinetic energy of the fluid disturbance, is highly concentrated at specific times.  
These spectacular displays of reversibility were captured experimentally for Vlasov, there known as \emph{plasma echoes}, in the work of \cite{MalmbergWharton68}. 
The analogous `Euler echoes' were recently studied and observed both numerically \cite{VMW98,Vanneste02} and experimentally \cite{YuDriscoll02,YuDriscollONeil}. 

The careful analysis of plasma echoes in the Vlasov equations is crucial in the proof of Mouhot and Villani \cite{MouhotVillani11}, as these are the dominant weakly nonlinear effect that could lead to instability.
Similarly, we also need to control this cascade for \eqref{def:2DEuler}, although our approach is quite different. 
To begin this analysis, first note that the basic challenge to the proof of Theorem \ref{thm:Main} is 
controlling the regularity of solutions to \eqref{Euler2}.
Since we must pay regularity to deduce decay on the velocity $u$, it is natural to consider the 
frequency interactions in the product $u\cdot \grad f$ with the frequencies of $u$ much larger than $f$.
This leads us to study a simpler model 
\begin{align} 
\partial_t f = -u\cdot \grad f_{lo}, \label{def:origlin}
\end{align}
where $f_{lo}$ is a given function that we think of as much smoother than $f$. 
As we see from \eqref{Euler2}, $u$ consists of several terms, however let us focus on the term we think should be the worst and also ignore the $v^\prime$, further reducing to the problem: 
\begin{align*} 
\partial_tf = \partial_v P_{\neq 0} \phi \partial_z f_{lo}.    
\end{align*} 
Suppose that instead of $f = \Delta_t\phi$, we had $f = \partial_{zz}\phi + (\partial_y - t\partial_z)^2\phi$ as in \eqref{orr-cri}, then the problem is actually linear and on the Fourier side:
\begin{align*} 
\partial_t\hat{f}(t,k,\eta) = \frac{1}{2\pi}\sum_{l \neq 0}\int_\xi \frac{\xi (k-l)}{l^2 + \abs{\xi - lt}^2}\hat{f}(l,\xi) \hat f_{lo}(t,k-l,\eta-\xi) d\xi.    
\end{align*}  
Since $f_{lo}$ weakens interactions between well-separated frequencies, let us consider a discrete model with $\eta$ as a fixed parameter:  
\begin{align} 
\partial_t\hat{f}(t,k,\eta) =  \frac{1}{2\pi}\sum_{l \neq 0} \frac{\eta (k-l)}{l^2 + \abs{\eta - lt}^2}\hat{f}(l,\eta)f_{lo}(t,k-l,0). 
\end{align}
As time advances this system of ODEs will go through resonances or ``critical times'' given by $t = \frac{\eta}{k}$, at which time the $k$ mode strongly forces the others.   
If $\abs{\eta}k^{-2} \ll 1$ then the critical time does not have a serious detriment; see \cite{BM13} for how this is dealt with. 
Henceforth only consider $\abs{\eta}k^{-2} > 1$.  
The scenario we are most concerned with is a high-to-low cascade in which the $k$ mode has a strong effect at time $\eta/k$ that excites the $k-1$ mode which has a strong effect at time $\eta/(k-1)$ that excites the $k-2$ mode and so on. 
Now focus near one critical time $\eta/k$ on a time interval of length roughly $\eta/k^2$, namely 
$ I_k= [ \eta/k -  \eta/k^2,  \eta/k +  \eta/k^2 ] $
 and consider the interaction between the mode $k$ and a nearby mode $l$ with $l \neq k$.
If one takes absolute values and retains only the leading order terms, then this reduces to the much simpler system of two ODEs (thinking of $f_{lo} = O(\kappa)$) which we refer to as the \emph{toy model}: 
\begin{subequations} \label{toy}
\begin{align}
\partial_tf_R & = \kappa\frac{k^2}{\abs{\eta}}f_{NR}, \\
\partial_tf_{NR} & = \kappa\frac{\abs{\eta}}{k^2 + \abs{\eta-kt}^2}f_{R}, 
\end{align}
\end{subequations} 
where we think of $f_R$ as being the evolution of the $k$ mode and $f_{NR}$ being the evolution of a nearby mode $l$ with $l \neq k$.
The factor $k^2/\abs{\eta}$ in the ODE for $f_R$ is an upper bound on the strongest interaction a non-resonant mode, for example the $k-1$ mode, can have with the resonant mode. 
Obviously \eqref{toy} represents a major simplification compared to \eqref{def:origlin}, however it will be sufficient to prove Theorem \ref{thm:Main}. 
 It is important to note that if at the beginning of the interval $I_k$, we have 
  $f_R = f_{NR}$, then   over the  interval  $I_k$, both 
$f_R$ and $f_{NR}$ are at most  amplified by roughly the same factor
$C  (\tfrac{\eta}{k^2})^{1+2C\kappa}$ (though they crucially are not amplified by the same amount 
on the left and right parts of the interval).  Taking the product of these amplifications 
for  $k = E(\sqrt{\eta}), E(\sqrt{\eta}) -1, ..., 1  $ yields a total amplification which is 
$ O(e^{C \sqrt{\eta}})$. 
This indicates that unless there is some special structure or cancellation not taken into account, 
the growth of high frequencies will cause a loss of Gevrey-2 regularity of the solution as $t \rightarrow \infty$. 
Therefore, in order to maintain control, the initial data must have at least this much regularity to lose, 
and this is the origin of the requirement $s > 1/2$ (or at least $s \geq 1/2$). 

The most important point of the toy model is not that it derives the expected regularity requirement (C. Mouhot and C. Villani have informed the authors that this same requirement can be derived in a manner similar to the approach used in \cite{MouhotVillani11}).  
Instead, the real use of the toy model will be in the design of our norm; see \S\ref{sec:MainEnergy} and \cite{BM13} for more information.

\subsection{Main Energy Estimates} \label{sec:MainEnergy}
The end goal is to prove that \eqref{Euler2} cannot behave significantly worse than that predicted by \eqref{toy}. 
It is insufficient to approximate the behavior of \eqref{toy} by an imprecise norm such as Gevrey-2, as this would result in either a rapid loss of all regularity control or a time growth much larger than the $O(t^{-2}\log t)$ damping in \eqref{Euler2} could hope to overpower.
Instead, we build \eqref{toy} \emph{into the energy estimate}. 
The key idea is the carefully designed time-dependent norm:
\begin{align*} 
\norm{A(t,\grad)f}^2_2 = \sum_k\int_\eta \abs{A_k(t,\eta)\hat f_k(t,\eta)}^2 d\eta,
\end{align*} 
where the multiplier $A$ has several components:
\begin{align*} 
A_k(t,\eta) = e^{\lambda(t)\abs{k,\eta}^s }  \jap{k,\eta}^\sigma  J_k(t,\eta)  B_k(t,\eta). 
\end{align*}
The goal is to get a \emph{uniform} bound on an energy like $\norm{A(t)f(t)}_2$; to do so we trade regularity for the necessary decay and indeed $A(t,k,\eta)$ becomes weaker and weaker as $t \rightarrow \infty$. 
The index $\lambda(t)$ is the `bulk' Gevrey$-\frac{1}{s}$ regularity and is chosen to satisfy 
\begin{subequations} \label{def:lambdat}
\begin{align} 
\lambda(t) & = \frac{3}{4}\lambda + \frac{1}{4}\lambda^\prime, \quad t \leq 1 \label{def:lambdashort} \\   
\dot\lambda(t) & = -K(\lambda^\prime,s,\sigma)\frac{\epsilon}{\jap{t}^{2\tilde q}}(1 + \lambda(t)), \quad t > 1
\end{align} 
\end{subequations}
for some $K(\lambda^\prime,s,\sigma)$ determined by the proof and $\tilde q = s/8 + 7/16$ (note $2\tilde q > 1$).
The reason for \eqref{def:lambdashort} is to account for the behavior of the solution for short time; see \cite{BM13} for this minor detail. 
Eventually, $\epsilon$ will be chosen sufficiently small such that 
 $\lambda(t) > (\lambda_0 + \lambda^\prime)/2$ for all $t>0$.
The use of a time-varying index of regularity is classical, for example the Cauchy-Kovalevskaya local existence theorem of Nirenberg \cite{Nirenberg72,Nishida77}; for more directly relevant works which use this kind of regularity loss, see \cite{FoiasTemam89,LevermoreOliver97,Chemin04,KukavicaVicol09,CGP11,MouhotVillani11}. 
The essential content of such time-varying regularity is that it allows high frequencies to grow faster than low frequencies.  
The Sobolev correction with $\sigma > 12$ fixed is included mostly for technical convenience.

The multiplier for dealing with the Orr mechanism and the associated nonlinear frequency cascade is
\begin{align} 
J_k(t,\eta) =   \frac{e^{\mu\abs{\eta}^{1/2}}}{w_k(t,\eta)}  + e^{\mu\abs{k}^{1/2}}, \label{def:J}
\end{align} 
where $ w_k(t,\eta)$ is constructed to approximate the behavior of \eqref{toy} near the critical times.
The key point is that $J_k(t,\eta)$ permits the solution to undergo growth consistent with \eqref{toy} (but not worse).  
An important consequence is that due to the behavior of \eqref{toy}, $J$ imposes \emph{more} regularity on modes which satisfy $t \sim \frac{\eta}{k}$ (the `resonant modes') than those that do not (the `non-resonant modes').
This discrepancy can create a useful gain or a dangerous loss of regularity when comparing ratios of the form $J_k(t,\eta)/J_l(t,\xi)$ (see \cite{BM13}). 

The multiplier $B_k(t,\eta)$ is given by the following, for $\gamma = s/2 + 1/4$,
\begin{align} 
B_k(t,\eta) = e^{\mu\abs{\eta}^{\gamma}}b(t,\eta)  + e^{\mu\abs{k}^{\gamma}}, \label{def:Bk}
\end{align} 
where $b(t,\eta)$ is defined in \cite{BM13} to account for an estimated loss of regularity which can occur from nonlinear interactions with the coefficients of $\Delta_t$.  
The second terms in $B_k(t,\eta)$ and $J_k(t,\eta)$ are to balance the regularities in $z$ and $v$, as a large discrepancy between the two would be inconsistent with the hyperbolic nature of the 2D Euler equations which makes no inherent structural distinction between the $x$ and $y$ directions (other than that imposed by the background flow).  

Eventually, our final aim is to prove the energy estimate
\begin{align} 
E(t) = \frac{1}{2}\norm{Af}_2^2 + \jap{t}^{4 - K_D\epsilon} \norm{\frac{A(t)}{\jap{\partial_v}^{s}}[\partial_t v]}_2^2 \lesssim \epsilon^2,  \label{ineq:energy}
\end{align}
where $K_D > 0$ is some constant fixed by the proof. 
Of course, there are really two coupled energy estimates: the one which controls $Af$ and the one which controls the convergence of the coordinate system $[\partial_t v]$.
The latter is performed using the momentum equation and contains a `hidden' structure which also seems to require $s \geq 1/2$ independently of the cascade implied by \eqref{toy}; see \cite{BM13} for more information on this interesting detail. 
As alluded to above, \eqref{ineq:energy} implies Theorem \ref{thm:Main} provided (say) $\lambda(t) > (\lambda +\lambda^\prime)/2$ and $\epsilon$ is chosen sufficiently small (see \cite{BM13}).  

A number of technical tools are employed to derive \eqref{ineq:energy}, which comprises the vast majority of the proof of Theorem \ref{thm:Main}. One tool certainly worth mentioning is para-differential calculus, specifically the paraproduct decomposition of \cite{Bony81}.
This is the natural Fourier-analytic tool for making some of the heuristics employed in the derivation of \eqref{toy} rigorous 
and provides a kind of linearization of high frequencies around lower frequencies.
In particular, our proof avoids the use of a Newton iteration scheme such as that employed in \cite{MouhotVillani11}, a technique which seems frustrated by the quasilinear aspect of \eqref{def:2DEuler}.
In particular, the paraproduct decomposes \eqref{Euler2} into three components: 
\begin{align} 
\partial_t f  + \mathcal T_{u}\grad f + \mathcal T_{\grad f} u + \mathcal R(u,\grad f) = 0,
\end{align}
where $\mathcal T_u \grad f$ denotes the part of the product $u\cdot \grad f$ where the frequencies of $f$ are much larger than that of $u$
and $\mathcal T_{\grad f} u$ denotes the part of the product $u\cdot \grad f$ where the frequencies of $u$ are much larger than that of $f$. 
We refer to the first term as `transport' since it preserves a transport-like structure and in analogy with \cite{MouhotVillani11}, we refer to the second term as `reaction'.   
Notice that the reaction term is exactly analogous to \eqref{def:origlin} in \S\ref{sec:Nonlinear}. 
Treating the transport contributions is done by adapting the classical methods for obtaining Gevrey regularity estimates for transport equations \cite{FoiasTemam89,LevermoreOliver97,KukavicaVicol09}.

\subsection{Elliptic estimates}
The high frequencies of the streamfunction $\phi$ appear in the reaction term in the previous section. 
With suitable smallness hypotheses, we do have an analogue of \eqref{ineq:phidecay} for $\Delta_t$ (in \cite{BM13}, we refer to it as the `lossy' elliptic estimate). 
However, in the reaction term we cannot apply this, since we would always need more regularity on $f$ than we know how to control (and as discussed, we cannot hope to improve \eqref{ineq:phidecay}). 
This unavoidable loss of regularity is reminiscent of the losses that occur in Nash-Moser iterations, but we will not be using an iterative approach. 

The `precision elliptic estimate' in \cite{BM13} is a primary technical component of our work (and the one with the most obscure physical interpretation).
The point is to invert the operator $\Delta_t$ and be able to gain some regularity/decay \emph{or} be able to precisely compensate the loss of ellipticity by weakening the norm in a corresponding fashion.  
Not coincidentally, due to how the loss of ellipticity arises in \eqref{toy}, this is exactly what $J$ does if instead of $\Delta_t$ we had the linear operator
\begin{align}
\Delta_L = \partial_{z}^2 + (\partial_v - t\partial_z)^2. \label{def:DeltaL}
\end{align} 
Instead, what happens is that by treating $\Delta_t$ as a perturbation of $\Delta_L$ we end up also with high frequency contributions from the coefficients $v^\prime$ and $v^{\prime\prime}$ in \eqref{def:Deltat}. 
These high frequency contributions become most dangerous where $\Delta_t$ loses ellipticity and it was precisely to compensate for these that the multiplier $B$ is included in the definition of $A$.

\section{Conclusion} 
The proof of Theorem \ref{thm:Main} is distinct from the proof of  Mouhot and Villani \cite{MouhotVillani11}, 
however see \cite{BM13} for a brief discussion of some of the important mathematical parallels and physical relationships.  
Although the Euler and Vlasov-Poisson systems have several fundamental differences,
 we are currently working with C. Mouhot  towards a simpler, alternative proof of the Landau damping result in \cite{MouhotVillani11}
for all Gevrey class smaller than three (e.g. $s > 1/3$) using some of the ideas of \cite{BM13}. 

 Orr and Kelvin (and many others) expressed doubt that the inviscid problem could be stable unless the set of permissible data was restricted in some way, predicting that the set of stable data would essentially vanish as the Reynolds number increased. 
While Theorem \ref{thm:Main} is in contradiction with this sentiment since it provides an open 
set (in Gevrey regularity) around zero for which the inviscid damping holds, 
their viewpoint could still be reconciled with our results. 
First, for inviscid flows we conjecture that in lower regularities than Gevrey-2, (one natural guess is $H^s$ for $s > 3/2$), `most' sufficiently small initial data will damp as in Theorem \ref{thm:Main} but there may exist a `non-generic' set of initial data which lead to unstable solutions that become trapped near cat's eye vortex-like solutions.
Second, we conjecture that for high Reynolds number flows, an appropriate analogue of Theorem \ref{thm:Main} still holds for  initial data $\omega^R_{in} + \omega_{in}^\nu$ where $\omega^R_{in}$ has Gevrey-$\frac{1}{s}$ regularity as in the Euler case 
 and $\omega_{in}^{\nu}$ has Sobolev regularity with 
 a small norm    which goes to zero when the Reynolds number goes to infinity.
 In particular this would put the result of \cite{BM13} as the inviscid limit of the high Reynolds number flow.
  
The viscous linearized problem was solved by Kelvin in 1887 \cite{Kelvin87} and is arguably the simplest example of `mixing-enhanced diffusion'. In this case, all flows return to Couette, so the analogue of Theorem \ref{thm:Main} would be a statement about the relative time-scales in the problem and how they behave in the inviscid limit. 
We are currently investigating the latter conjecture (the former appears far out of reach for now).  
    
Mixing is a fundamental mechanism in fluid mechanics and plasma physics that 
permits conservative, reversible systems to exhibit decay due to mixing in phase space.
There are many settings in which 
the linearized problem  predicts asymptotic stability by mixing, however each seems to require non-trivial new ideas to approach on the nonlinear level.

The most obvious related problem is the inviscid damping of a 
 general class of shear flows for which the linearization predicts asymptotic stability.
However, more general shears fundamentally change the structure of the critical times used in the case of the Couette flow.  
 The associated nonlinear effects are more complicated to control and 
 would require some precise adaptations to handle.  
For example, we currently believe that much of the Fourier analysis would need to be replaced by more delicate pseudo-differential calculus. 

A related extension would be to study the problem in the presence of no-penetration boundaries in 
 the vertical direction $y$, which may shed light onto the way boundaries could cause some of the  instabilities seen in experiments.     
A third problem would be to remove the periodicity in $x$, altering the physical mechanism from mixing to \emph{filamentation}. In this case, the decay may be caused by a combination of filamentation and  dispersion. 

Other examples in fluid mechanics which involve stability by mixing include the $\beta$-plane model  \cite{Boyd83,Tung83}, stratified shear flows \cite{Majda03,CV13} and the particularly fundamental `vortex axisymmetrization' problem. 
See \cite{Gilbert88,GilbertBassom98,SchecterEtAl00,BMT13,YuDriscoll02,YuDriscollONeil} for a small piece of the extensive literature
 on this problem which dates back even to Rayleigh \cite{Rayleigh80} and Orr \cite{Orr07}. 
We hope that the methods proposed in our work are flexible and accessible enough to be extended to solve some of these other fundamental problems and shed further light on the important physical mechanisms behind them.

\bibliographystyle{plain} \bibliography{eulereqns}

\begin{thebibliography}{10}

\bibitem{BaggettEtAl}
J.S. Baggett, T.A. Driscoll, and L.N. Trefethen.
\newblock A mostly linear model of transition of turbulence.
\newblock {\em Phys. Fluids}, 7:833--838, 1995.

\bibitem{BM02}
N.~J. Balmforth and P.~J. Morrison.
\newblock Hamiltonian description of shear flow.
\newblock In {\em Large-scale atmosphere-ocean dynamics, {V}ol.\ {II}}, pages
  117--142. Cambridge Univ. Press, Cambridge, 2002.

\bibitem{BM95}
N.J. Balmforth and P.J. Morrison.
\newblock Normal modes and continuous spectraa.
\newblock {\em Annals of the New York Academy of Sciences}, 773(1):80--94,
  1995.

\bibitem{BMSEI95}
N.J. Balmforth and P.J. Morrison.
\newblock Singular eigenfunctions for shearing fluids {I}.
\newblock {\em Institute for Fusion Studies Report, University of
  Texas-Austin}, (692):1--80, 1995.

\bibitem{BMT13}
N.J. Balmforth, P.J. Morrison, and J.-L. Thiffeault.
\newblock Pattern formation in {Hamiltonian} systems with continuous spectra; a
  normal-form single-wave model.
\newblock {\em preprint}, 2013.

\bibitem{GilbertBassom98}
Andrew~P. Bassom and Andrew~D. Gilbert.
\newblock The spiral wind-up of vorticity in an inviscid planar vortex.
\newblock {\em J. Fluid Mech.}, 371:109--140, 1998.

\bibitem{BM13}
J.~Bedrossian and N.~Masmoudi.
\newblock Inviscid damping and the asymptotic stability of planar shear flows
  in the {2D Euler} equations.
\newblock {\em {arXiv:1306.5028}}.

\bibitem{Bony81}
J.M. Bony.
\newblock Calcul symbolique et propagation des singularit\'es pour les
  \'equations aux d\'eriv\'ees partielles non lin\'aires.
\newblock {\em Ann.Sc.E.N.S.}, 14:209--246, 1981.

\bibitem{BouchetMorita10}
F.~Bouchet and H.~Morita.
\newblock Large time behavior and asymptotic stability of the {2D Euler} and
  linearized {Euler} equations.
\newblock {\em Physica D}, 239:948--966, 2010.

\bibitem{Boyd83}
John~P Boyd.
\newblock The continuous spectrum of linear couette flow with the beta effect.
\newblock {\em Journal of the atmospheric sciences}, 40(9):2304--2308, 1983.

\bibitem{Briggs70}
R.J. Briggs, J.D. Daugherty, and R.H. Levy.
\newblock Role of {Landau} damping in crossed-field electron beams and inviscid
  shear flow.
\newblock {\em Phys. Fl.}, 13(2), 1970.

\bibitem{CagliotiMaffei98}
E.~Caglioti and C.~Maffei.
\newblock Time asymptotics for solutions of {Vlasov-Poisson} equation in a
  circle.
\newblock {\em J. Stat. Phys.}, 92(1/2), 1998.

\bibitem{CV13}
R~Camassa and C~Viotti.
\newblock Transient dynamics by continuous-spectrum perturbations in stratified
  shear flows.
\newblock {\em Journal of Fluid Mechanics}, 717, 2013.

\bibitem{Chemin04}
Jean-Yves Chemin.
\newblock Le syst\`eme de {N}avier-{S}tokes incompressible soixante dix ans
  apr\`es {J}ean {L}eray.
\newblock In {\em Actes des {J}ourn\'ees {M}ath\'ematiques \`a la {M}\'emoire
  de {J}ean {L}eray}, volume~9 of {\em S\'emin. Congr.}, pages 99--123. Soc.
  Math. France, Paris, 2004.

\bibitem{CGP11}
Jean-Yves Chemin, Isabelle Gallagher, and Marius Paicu.
\newblock Global regularity for some classes of large solutions to the
  {N}avier-{S}tokes equations.
\newblock {\em Ann. of Math. (2)}, 173(2):983--1012, 2011.

\bibitem{ConstantinEtAl08}
P.~Constantin, A.~Kiselev, L.~Ryzhik, and A.~Zlato{\v{s}}.
\newblock Diffusion and mixing in fluid flow.
\newblock {\em Ann. of Math. (2)}, 168(2):643--674, 2008.

\bibitem{cycon1987}
HL~Cycon, RG~Froese, W~Kirsch, and B~Simon.
\newblock Schr{\'o}dinger operators.
\newblock {\em Lecture Notes in Physics, Berlin-Heidelberg-New york: Springer
  Verlag}, 1987.

\bibitem{Degond86}
P.~Degond.
\newblock Spectral theory of the linearized {Vlasov-Poisson} equation.
\newblock {\em Trans. Amer. Math. Soc.}, 294(2):435--453, 1986.

\bibitem{FoiasTemam89}
C.~Foias and R.~Temam.
\newblock Gevrey class regularity for solutions of the {Navier-Stokes}
  equations.
\newblock {\em J. Funct. Anal.}, 87:359--369, 1989.

\bibitem{GMS12}
P.~Germain, N.~Masmoudi, and J.~Shatah.
\newblock Global solutions for the gravity water waves equation in dimension 3.
\newblock {\em Ann. of Math. (2)}, 175(2):691--754, 2012.

\bibitem{Gevrey18}
Maurice Gevrey.
\newblock Sur la nature analytique des solutions des \'equations aux
  d\'eriv\'ees partielles. {P}remier m\'emoire.
\newblock {\em Ann. Sci. \'Ecole Norm. Sup. (3)}, 35:129--190, 1918.

\bibitem{Gilbert88}
Andrew~D. Gilbert.
\newblock Spiral structures and spectra in two-dimensional turbulence.
\newblock {\em J. Fluid Mech.}, 193:475--497, 1988.

\bibitem{GV00}
J.~Ginibre and G.~Velo.
\newblock Long range scattering and modified wave operators for some {H}artree
  type equations. {I}.
\newblock {\em Rev. Math. Phys.}, 12(3):361--429, 2000.

\bibitem{GSV13}
Nathan Glatt-Holtz, Vladimir Sverak, and Vlad Vicol.
\newblock On inviscid limits for the stochastic navier-stokes equations and
  related models.
\newblock {\em arXiv preprint arXiv:1302.0542}, 2013.

\bibitem{HagstromMorrison}
G.I. Hagstrom and P.J. Morrison.
\newblock {Caldeira-Leggett} model, {Landau} damping and the {Vlasov-Poisson}
  system.
\newblock {\em Phys. D.}, 240:1652--1660, 2011.

\bibitem{HwangVelazquez09}
H.~J. Hwang and J.~J.~L. Vela{\'z}quez.
\newblock On the existence of exponentially decreasing solutions of the
  nonlinear {Landau} damping problem.
\newblock {\em Indiana Univ. Math. J}, pages 2623--2660, 2009.

\bibitem{IonescuPusateri13}
A.~Ionescu and F.~Pusateri.
\newblock Global solutions for the gravity water waves in {2D}.
\newblock {\em {arXiv:1303.5357}}.

\bibitem{Kelvin87}
Lord Kelvin.
\newblock Stability of fluid motion-rectilinear motion of viscous fluid between
  two parallel plates.
\newblock {\em Phil. Mag.}, (24):188, 1887.

\bibitem{Kraichnan67}
R.H. Kraichnan.
\newblock Inertial ranges in two-dimensional turbulence.
\newblock {\em Phys. Fluids}, 10(7), 1967.

\bibitem{KukavicaVicol09}
Igor Kukavica and Vlad Vicol.
\newblock On the radius of analyticity of solutions to the three-dimensional
  {E}uler equations.
\newblock {\em Proc. Amer. Math. Soc.}, 137(2):669--677, 2009.

\bibitem{Landau46}
L.~Landau.
\newblock On the vibration of the electronic plasma.
\newblock {\em J. Phys. USSR}, 10(25), 1946.

\bibitem{LevermoreOliver97}
D.~Levermore and M.~Oliver.
\newblock Analyticity of solutions for a generalized {Euler} equation.
\newblock {\em J. Diff. Eqns.}, 133:321--339, 1997.

\bibitem{LiLin11}
Y.~Charles Li and Zhiwu Lin.
\newblock A resolution of the {S}ommerfeld paradox.
\newblock {\em SIAM J. Math. Anal.}, 43(4):1923--1954, 2011.

\bibitem{LinZeng11}
Z.~Lin and C.~Zeng.
\newblock Inviscid dynamic structures near {Couette} flow.
\newblock {\em Arch. Rat. Mech. Anal.}, 200:1075--1097, 2011.

\bibitem{LZ11b}
Zhiwu Lin and Chongchun Zeng.
\newblock Small {BGK} waves and nonlinear {L}andau damping.
\newblock {\em Comm. Math. Phys.}, 306(2):291--331, 2011.

\bibitem{Lindzen88}
R.~Lindzen.
\newblock Instability of plane parallel shear flow (toward a mechanistic
  picture of how it works).
\newblock {\em PAGEOPH}, 126(1), 1988.

\bibitem{Majda03}
Andrew Majda.
\newblock {\em Introduction to {PDE}s and waves for the atmosphere and ocean},
  volume~9 of {\em Courant Lecture Notes in Mathematics}.
\newblock New York University Courant Institute of Mathematical Sciences, New
  York, 2003.

\bibitem{MalmbergWharton68}
J.~Malmberg, C.~Wharton, C.~Gould, and T.~O'Neil.
\newblock Plasma wave echo.
\newblock {\em Phys. Rev. Lett.}, 20(3):95--97, 1968.

\bibitem{Morrison98}
P.~J. Morrison.
\newblock Hamiltonian description of the ideal fluid.
\newblock {\em Rev. Modern Phys.}, 70(2):467--521, 1998.

\bibitem{Morrison00}
P.~J. Morrison.
\newblock Hamiltonian description of vlasov dynamics: action-angle variables
  for the continuous spectrum.
\newblock {\em Trans. Theory and Stat. Phys.}, 29(3-5):397--414, 2000.

\bibitem{MouhotVillani11}
C.~Mouhot and C.~Villani.
\newblock On {Landau} damping.
\newblock {\em Acta Math.}, 207:29--201, 2011.

\bibitem{Nakanishi02}
Kenji Nakanishi.
\newblock Modified wave operators for the {H}artree equation with data, image
  and convergence in the same space.
\newblock {\em Commun. Pure Appl. Anal.}, 1(2):237--252, 2002.

\bibitem{Nirenberg72}
L.~Nirenberg.
\newblock An abstract form of the nonlinear {Cauchy-Kowalewski} theorem.
\newblock {\em J. Diff. Geom.}, 6:561--576, 1972.

\bibitem{Nishida77}
T.~Nishida.
\newblock A note on a theorem of {Nirenberg}.
\newblock {\em J. Diff. Geom.}, 12:629--633, 1977.

\bibitem{Orr07}
W.~Orr.
\newblock The stability or instability of steady motions of a perfect liquid
  and of a viscous liquid, {Part I}: a perfect liquid.
\newblock {\em Proc. Royal Irish Acad. Sec. A: Math. Phys. Sci.}, 27:9--68,
  1907.

\bibitem{Rayleigh80}
Lord Rayleigh.
\newblock On the {S}tability, or {I}nstability, of certain {F}luid {M}otions.
\newblock {\em Proc. London Math. Soc.}, S1-11(1):57, 1880.

\bibitem{ReedSimonIII}
M.~Reed and B.~Simon.
\newblock {\em Methods of Modern Mathematical Physics {III}: Scattering
  theory}.
\newblock Academic Press, 1979.

\bibitem{SchecterEtAl00}
D.A. Schecter, D.~Dubin, A.C. Cass, C.F. Driscoll, and I.M.~Lansky et. al.
\newblock Inviscid damping of asymmetries on a two-dimensional vortex.
\newblock {\em Phys. Fl.}, 12, 2000.

\bibitem{Shnirelman12}
Alexander Shnirelman.
\newblock On the long time behavior of fluid flows.
\newblock {\em preprint}, 2012.

\bibitem{StrogatzEtAl92}
S.H. Strogatz, R.E. Mirollow, and P.C. Matthews.
\newblock Coupled nonlinear oscillators below the synchronization threshold:
  relaxation by generalized {Landau} damping.
\newblock {\em Phys. Rev. Let.}, 68(18):2730--2733, 1992.

\bibitem{TataronisGrossman73}
J.~Tataronis and W.~Grossmann.
\newblock Decay of {MHD} waves by phase mixing.
\newblock {\em Z. Physik}, 261:203--216, 1973.

\bibitem{TTRD93}
Lloyd~N. Trefethen, Anne~E. Trefethen, Satish~C. Reddy, and Tobin~A. Driscoll.
\newblock Hydrodynamic stability without eigenvalues.
\newblock {\em Science}, 261(5121):578--584, 1993.

\bibitem{Tung83}
K.K. Tung.
\newblock Initial-value problems for {Rossby} waves in a shear flow with
  critical level.
\newblock {\em J. Fluid Mech.}, 133:443--469, 1983.

\bibitem{VKampen55}
N.G. van Kampen.
\newblock On the theory of stationary waves in plasmas.
\newblock {\em Physica}, 21:949--963, 1955.

\bibitem{Vanneste02}
J.~Vanneste.
\newblock Nonlinear dynamics of anisotropic disturbances in plane {C}ouette
  flow.
\newblock {\em SIAM J. Appl. Math.}, 62(3):924--944 (electronic), 2002.

\bibitem{VMW98}
J~Vanneste, P.J. Morrison, and T~Warn.
\newblock Strong echo effect and nonlinear transient growth in shear flows.
\newblock {\em Physics of Fluids}, 10:1398, 1998.

\bibitem{YuDriscoll02}
J.H. Yu and C.F. Driscoll.
\newblock Diocotron wave echoes in a pure electron plasma.
\newblock {\em {IEEE} Trans. Plasma Sci.}, 30(1), 2002.

\bibitem{YuDriscollONeil}
J.H. Yu, C.F. Driscoll, and T.M. O`Neil.
\newblock Phase mixing and echoes in a pure electron plasma.
\newblock {\em Phys. of Plasmas}, 12(055701), 2005.

\end{thebibliography}

\end{document}